\documentclass[titlepage,draft,12pt]{article} 
\usepackage{amsfonts,latexsym,pstricks,pst-plot} 
\usepackage[a4paper]{geometry}

\newcommand{\re}{{\mathbb{R}}}

\newcommand{\uep}{u_{\ep}}

\newcommand{\wep}{w_{\ep}}
\newcommand{\Dep}{D_{\ep}}
\newcommand{\DDep}{\widehat{D}_{\ep}}

\newcommand{\Fep}{F_{\ep}}
\newcommand{\Gep}{G_{\ep}}
\newcommand{\Hep}{H_{\ep}}
\newcommand{\Pep}{P_{\ep}}
\newcommand{\Qep}{Q_{\ep}}
\newcommand{\Rep}{R_{\ep}}

\newcommand{\auq}[1]{|A^{1/2}#1|^{2}}

\newcommand{\dau}{D(A^{1/2})}
\newcommand{\au}[1]{|A^{1/2}#1|}
\newcommand{\da}{D(A)}
\newcommand{\ep}{\varepsilon}
\newcommand{\prf}{{\sc Proof.}$\;$}
\newcommand{\qed}{{\penalty 10000\mbox{$\quad\Box$}}}

\newcommand{\m}[1]{m(\au{#1}^{2})}

\newtheorem{thm}{Theorem}[section]

\newtheorem{prop}[thm]{Proposition}

\newtheorem{lemma}[thm]{Lemma}

\newtheorem{rmk}[thm]{Remark}

\title{Mildly degenerate Kirchhoff equations with weak dissipation:
global existence and time decay}
\author{Marina Ghisi\vspace{1ex}\\ {\normalsize
Universit\`a degli Studi di Pisa} \\{\normalsize Dipartimento di
Matematica ``Leonida Tonelli''}\\
{\normalsize 
PISA (Italy)}\\  
{\normalsize e-mail: \texttt{ghisi@dm.unipi.it}}\and 
Massimo Gobbino\vspace{1ex}\\ {\normalsize Universit\`a degli Studi di Pisa} 
\\{\normalsize Dipartimento di Matematica Applicata ``Ulisse Dini''}\\ 
{\normalsize 
 PISA (Italy)}\\  
{\normalsize e-mail: \texttt{m.gobbino@dma.unipi.it}}}
\date{}

\begin{document}
\maketitle
\begin{abstract}
	We consider the hyperbolic-parabolic singular perturbation problem
	for a \emph{degenerate} quasilinear Kirchhoff equation with
	\emph{weak} dissipation. This means that the coefficient of the
	dissipative term tends to zero when $t\to +\infty$. 
	
	We prove that the hyperbolic problem has a unique global solution
	for suitable values of the parameters.  We also prove that the
	solution decays to zero, as $t\to +\infty$, with the same rate of
	the solution of the limit problem of parabolic type.
	
\vspace{1cm}

\noindent{\bf Mathematics Subject Classification 2000 (MSC2000):}
35B25, 35B40, 35L70, 35L80.

\vspace{1cm} \noindent{\bf Key words:} singular perturbation,
degenerate Kirchhoff equations, quasilinear hyperbolic equation, weak
dissipation, energy estimates.
\end{abstract}
 
\section{Introduction}
Let $H$ be a real Hilbert space.  For every $x$ and $y$ in $H$, $|x|$
denotes the norm of $x$, and $\langle x,y\rangle$ denotes the scalar
product of $x$ and $y$.  Let $A$ be a self-adjoint linear operator on
$H$ with dense domain $D(A)$.  We assume that $A$ is nonnegative,
namely $\langle Ax,x\rangle\geq 0$ for every $x\in D(A)$, so that for
every $\alpha\geq 0$ the power $A^{\alpha}x$ is defined provided that
$x$ lies in a suitable domain $D(A^{\alpha})$.

We consider the Cauchy problem 
\begin{equation}
	\ep\uep''(t)+\frac{1}{(1+t)^{p}}\uep'(t)+
	\au{\uep(t)}^{2\gamma}A\uep(t)=0 \hspace{2em}
	\forall t\geq 0,
	\label{pbm:h-eq}
\end{equation}
\begin{equation}
	\uep(0)=u_0,\hspace{2em}\uep'(0)=u_1,
	\label{pbm:h-data}
\end{equation}
where $\ep>0$, $p\geq 0$, and $\gamma> 0$.
Equation (\ref{pbm:h-eq}) is the prototype of all degenerate Kirchhoff
equations with weak dissipation
\begin{equation}
	\ep\uep''(t)+b(t)\uep'(t)+
	\m{\uep(t)}A\uep(t)=0 \hspace{2em}
	\forall t\geq 0,
	\label{pbm:h-eq-gen}
\end{equation}
where $b:[0,+\infty)\to(0,+\infty)$ and $m:[0,+\infty)\to [0,+\infty)$
are given functions which are always assumed to be of class $C^{1}$ (or
at least locally Lipschitz continuous), unless otherwise stated.  It is
well known that (\ref{pbm:h-eq}) is the abstract setting of a
quasilinear nonlocal partial differential equation of hyperbolic type
which was proposed as a model for small vibrations of strings and
membranes.

Equation (\ref{pbm:h-eq-gen}) is called nondegenerate (or strictly
hyperbolic) when
$$\mu:=\inf_{\sigma\geq 0}m(\sigma)>0,$$
and mildly degenerate when $\mu=0$ but $\m{u_{0}}\neq 0$. In the
special case of equation (\ref{pbm:h-eq}) this assumption reduces to 
\begin{equation}
	A^{1/2}u_{0}\neq 0.
	\label{hp:mdg}
\end{equation}

Concerning the dissipation term $b(t)\uep'(t)$, we have constant
dissipation when $b(t)\equiv\delta>0$ is a positive constant, and weak
dissipation when $b(t)\to 0$ as $t\to +\infty$.  Finally, the operator
$A$ is called coercive when
\begin{equation}
	\nu:=\inf\left\{\frac{\langle Ax,x\rangle}{|x|^{2}}:x\in D(A),\
	x\neq 0\right\}>0,
    \label{hp:coercive}
\end{equation}
and noncoercive when $\nu=0$.

The singular perturbation problem in its generality consists in
proving the convergence of solutions of (\ref{pbm:h-eq-gen}),
(\ref{pbm:h-data}) to solutions of the first order problem
\begin{equation}
	b(t)\uep'(t)+ \m{\uep(t)}A\uep(t)=0, \hspace{2em}
	u(0)=u_{0},
	\label{pbm:par}
\end{equation}
obtained setting formally $\ep=0$ in (\ref{pbm:h-eq-gen}), and
omitting the second initial condition in~(\ref{pbm:h-data}).

The singular perturbation problem gives rise to several subproblems.
The first step is of course the existence of global solutions for the
limit problem (\ref{pbm:par}).  This has been established in
\cite{k-par} under very general assumptions.  The second step is the
existence of a global solution for the hyperbolic problem
(\ref{pbm:h-eq-gen}), (\ref{pbm:h-data}).  The third step is the
convergence of solutions $\uep(t)$ of the hyperbolic problem to the
solution $u(t)$ of the parabolic problem.  The final goal are the so
called error-decay estimates which prove in the same time that the
difference $\uep(t)-u(t)$ decays to 0 as $t\to +\infty$ (with the same
rate of $u(t)$), and tends to 0 as $\ep\to 0^{+}$.


The second and third step have generated a considerable literature,
which we sum up below.

\subparagraph{\emph{\textmd{Nondegenerate Kirchhoff equations with
constant dissipation}}}

The first results where obtained in the eighties by \textsc{E.\ De
Brito~\cite{debrito}} and \textsc{Y.\ Yamada}~\cite{yamada}.  They
independently proved the global solvability of the hyperbolic problem
with initial data $(u_{0},u_{1})\in D(A)\times D(A^{1/2})$ under a
suitable assumption involving $\ep$, the initial data, and the
constant dissipation $\delta$.  Once that $\delta$ and the initial
data are fixed, this condition holds true provided that $\ep$ is small
enough.  The key step in their proofs, as well as in all the
subsequent literature, is to show that solutions satisfy an a priori
estimate such as
\begin{equation}
	\ep\: \frac{\left|m'(\auq{\uep(t)})\right|}{\m{\uep(t)}}
	\cdot|\uep'(t)|\cdot|A\uep(t)|\leq b(t),
	\label{a-priori}
\end{equation}
which is clearly more likely to be true when $\ep$ is small enough.
Existence of global solutions without the smallness assumption on
$\ep$ remains a challenging open problem, as well as the
nondissipative case $b(t)\equiv 0$.

More recently \textsc{H.\ Hashimoto} and \textsc{T.\
Yamazaki}~\cite{yamazaki} obtained optimal error-decay estimates for
the singular perturbation problem.  Thanks to these estimates, which
improve or extend all previous works (see \cite{ew,k-cattaneo}), this
case can be considered quite well understood.

\subparagraph{\emph{\textmd{Degenerate Kirchhoff equations with constant
dissipation}}}

The case where $m(\sigma)=\sigma^{\gamma}$ (with $\gamma\geq 1$) has
been studied in the nineties by \textsc{K.\ Nishihara} and
\textsc{Y.\ Yamada}~\cite{ny} (see also \cite{ono-mm}).  The result is
the existence of a unique global solution for the mildly degenerate
equation provided that $\ep$ is small enough.  Later on this existence
result was extended by the authors \cite{gg:k-dissipative} to
arbitrary locally Lipschitz continuous nonlinearities $m(\sigma)\geq
0$, and by the first author \cite{ghisi1, g:non-lip} to non-Lipschitz
nonlinearities of the form $m(\sigma)=\sigma^{\gamma}$ with
$\gamma\in(0,1)$.

All the quoted papers considered also the asymptotic behavior of
solutions, but the estimates proved therein where in general far from
being optimal.  In the meanwhile sharp decay estimates were the
subject of a series of papers by \textsc{T.\ Mizumachi}
(\cite{mizu-ade,mizu-nc}) and \textsc{K.\ Ono}
(\cite{ono-kyushu,ono-aa}), in which however only the special case
$m(\sigma)=\sigma$ was considered.  More recently the
authors~\cite{gg:k-decay} obtained optimal and $\ep$-independent decay
estimates for the general case.  As expected the result is that
solutions of the hyperbolic problem always decay as the corresponding
solutions of the limit problem.

As for the singular perturbation problem, in~\cite{gg:k-PS} the
authors proved that $\uep(t)\to u(t)$ uniformly in time, but without
sharp error-decay estimates, which in this case remain an open
problem.

From the technical point of view, the difficulty is that in the
degenerate case the denominator in (\ref{a-priori}) may vanish.  This
cannot happen for $t=0$ due to the mild nondegeneracy assumption, but
it does happen in the limit as $t\to +\infty$ due to the decay of 
solutions.

\subparagraph{\emph{\textmd{Nondegenerate Kirchhoff equations with
weak dissipation}}}

Let us come to the nondegenerate case with $b(t)=(1+t)^{-p}$.  What
complicates things is the competition between the smallness of $\ep$
and the smallness of $b(t)$.  In particular it is no more enough to
prove that the left-hand side of (\ref{a-priori}) is bounded, but it
is necessary to prove that it decays with an a-priori fixed rate.

This is the reason why this problem has been solved only in recent
years in some papers by {\sc M.\ Nakao} and {\sc J.\
Bae}~\cite{nakao}, by \textsc{T.\
Yamazaki}~\cite{yamazaki-wd,yamazaki-cwd}, and by the
authors~\cite{gg:cwd}.  The result is that for every $p\in[0,1]$, and
every $(u_{0},u_{1})\in\da\times\dau$, the hyperbolic problem has a
unique global solution provided that $\ep$ is small enough.  Moreover
the solution decays to $0$ as $t\to +\infty$ as the solution of the
limit problem, and optimal error-decay estimates for the singular
perturbation can be proved.  When $p>1$ existence of global solutions
for the hyperbolic problem is still an open problem, but in any case
solutions cannot decay to $0$ as $t\to +\infty$.  On the other hand,
solutions of the limit parabolic problem decay to zero also for $p>1$,
faster and faster as $p$ grows.

This means that a threshold appears.  When $p\in[0,1]$ the smallness
of $\ep$ is dominant over the smallness of $b(t)$, and
(\ref{pbm:h-eq-gen}) behaves like a parabolic equation.  When $p>1$
the smallness of $b(t)$ is dominant over the smallness of $\ep$, and
(\ref{pbm:h-eq-gen}) behaves like a nondissipative hyperbolic
equation.

\subparagraph{\emph{\textmd{Degenerate Kirchhoff equations with weak
dissipation}}}

Let us finally come to equation (\ref{pbm:h-eq}), which is the object
of this paper.  Now in (\ref{a-priori}) the smallness of $\ep$ has to
compete both with the decay of $b(t)$, and with the vanishing of the
denominator.  So one needs a priori decay estimates for terms whose
denominator vanishes in the limit.

To our knowledge the only previous result for this equation was
obtained at the end of the nineties by \textsc{K.\
Ono}~\cite{ono-wd}, who proved global existence for the mildly
degenerate case when $\gamma=1$, $p\in[0,1/3]$, and of course $\ep$ is
small enough.  We recall that $m(\sigma)=\sigma$ is the only case
where sharp decay estimates were already available in those years.

In this paper we consider the global solvability of (\ref{pbm:h-eq}),
(\ref{pbm:h-data}) with more general values of the parameters, and
initial data $(u_{0},u_{1})\in\da\times\dau$ satisfying
(\ref{hp:mdg}).  

Our first result (Theorem~\ref{thm:c-existence}) concerns the coercive
case.  Under this assumption we prove that a unique global solution
exists provided that $\gamma>0$, $p\in[0,1]$ and $\ep$ is small
enough.  We also prove sharp decay estimates as $t\to +\infty$.

Our second result (Theorem~\ref{thm:nc-existence}) concerns the
noncoercive case.  In this case we prove that a unique global solution
exists provided that $\gamma\geq 1$,
$p\in[0,(\gamma^{2}+1)/(\gamma^{2}+2\gamma-1)]$, and $\ep$ is small
enough.  Note that the supremum of this interval is 1 both when
$\gamma=1$ and when $\gamma\to +\infty$, but it is strictly smaller
than 1 for $\gamma>1$.  We also provide decay estimates for solutions.

Finally in both cases we show (Theorem~\ref{thm:p>1}) that for $p>1$
solutions of (\ref{pbm:h-eq}) (provided that they exist, which remains
an open problem) do not decay to 0 as $t\to +\infty$.

From the point of view of global solvability and decay properties
these results show that in the coercive case equation (\ref{pbm:h-eq})
behaves like the nondegenerate one, exhibiting nondissipative
hyperbolic behavior for $p>1$, and parabolic behavior for $p\in[0,1]$.
We point out that this is true also in the non-Lipschitz case
$\gamma\in(0,1)$.  In the noncoercive case (with $\gamma\geq 1$) we
have once again hyperbolic behavior for $p>1$, and parabolic behavior
for $p\in[0,(\gamma^{2}+1)/(\gamma^{2}+2\gamma-1)]$. 

Proofs rely on the techniques introduced in~\cite{gg:k-decay} in order
to prove sharp decay estimates.  When the operator is coercive the
decay rate depends only on $p$ and $\gamma$.  When the operator is
noncoercive the decay rate belongs to a range depending on $p$ and
$\gamma$, but within this range it seems to depend on the initial
conditions.  The existence of a range of possible decay rates is what
in the noncoercive case creates the no-man's land between
$(\gamma^{2}+1)/(\gamma^{2}+2\gamma-1)$ and 1.  What happens when $p$
lies in this interval is not clear yet.

In this paper we don't consider the behavior of solutions as $\ep\to
0^{+}$, even if all our $\ep$-independent estimates are for sure a
first step in this direction.  We just mention that a simple
adaptation of the arguments of \cite{k-cattaneo} and \cite{gg:k-PS}
should be enough to prove two types of result: that $\uep\to u$
uniformly in time (without estimates of the convergence rate), and
that $\uep\to u$ in every interval $[0,T]$ with an estimate of the
error depending on $T$.  On the other hand, obtaining error-decay
estimates analogous to the nondegenerate case seems to be a much more
difficult task.  Apart from the partial results of~\cite{gg:k-PS} this
problem is still open in the degenerate case, both with constant and
with weak dissipation.

\setcounter{equation}{0}
\section{Statements}\label{sec:statements}

Our first result concerns the global solvability of the hyperbolic
problem and decay properties of solutions in the case of coercive
operators.

\begin{thm}\label{thm:c-existence}
	Let $H$ be a Hilbert space, and let $A$ be a nonnegative
	selfadjoint (unbounded) operator with dense domain.  Let us assume
	that $A$ satisfies the coerciveness condition (\ref{hp:coercive}).
	Let $\gamma>0$, and let $p\in[0,1]$.  Let us assume that
	$(u_{0},u_{1})\in\da\times\dau$ satisfy (\ref{hp:mdg}).
	
	Then there exists $\ep_{0}>0$ such that for every
	$\ep\in(0,\ep_{0})$ problem (\ref{pbm:h-eq}), (\ref{pbm:h-data})
	has a unique global solution 
	\begin{equation}
		\uep\in C^{2}([0,+\infty);H)\cap C^{1}([0,+\infty);\dau) \cap
		C^{0}([0,+\infty);\da).
		\label{th:uep-reg}
	\end{equation}
	
	Moreover there exist positive constants $C_{1}$ and $C_{2}$ such that
	\begin{equation}
		\frac{C_{1}}{(1+t)^{(p+1)/\gamma}}\leq
		\auq{\uep(t)}\leq
		\frac{C_{2}}{(1+t)^{(p+1)/\gamma}}
		\quad\quad\forall t\geq 0;
		\label{th:decay-a1/2u}
	\end{equation}
	\begin{equation}
		\frac{C_{1}}{(1+t)^{(p+1)/\gamma}}\leq
		|A\uep(t)|^{2}\leq
		\frac{C_{2}}{(1+t)^{(p+1)/\gamma}}
		\quad\quad\forall t\geq 0;
		\label{th:decay-au}
	\end{equation}
	\begin{equation}
		|\uep'(t)|^{2}\leq
		\frac{C_{2}}{(1+t)^{2+(p+1)/\gamma}}
		\quad\quad\forall t\geq 0.
		\label{th:decay-u'}
	\end{equation}
\end{thm}

Our second result in the counterpart of Theorem~\ref{thm:c-existence}
in the case of noncoercive operators.

\begin{thm}\label{thm:nc-existence}
	Let $H$ be a Hilbert space, and let $A$ be a nonnegative
	selfadjoint (unbounded) operator with dense domain.  Let
	$\gamma\geq 1$, and let
	\begin{equation}
		0\leq p\leq \frac{\gamma^{2}+1}{\gamma^{2}+2\gamma-1}.  
		\label{hp:pnc}
	\end{equation}
	
	Let us assume
	that $(u_{0},u_{1})\in\da\times\dau$ satisfy (\ref{hp:mdg}).
	
	Then there exists $\ep_{0}>0$ such that for every
	$\ep\in(0,\ep_{0})$ problem (\ref{pbm:h-eq}), (\ref{pbm:h-data})
	has a unique global solution satisfying (\ref{th:uep-reg}).
	
	Moreover there exist constants $C_{1}$ and $C_{2}$ such that
	\begin{equation}
		\frac{C_{1}}{(1+t)^{(p+1)/\gamma}}\leq
		\auq{\uep(t)}\leq
		\frac{C_{2}}{(1+t)^{(p+1)/(\gamma+1)}}
		\quad\quad\forall t\geq 0;
		\label{th:decay-a1/2u-nc}
	\end{equation}
	\begin{equation}
		|A\uep(t)|^{2}\leq
		\frac{C_{2}}{(1+t)^{(p+1)/\gamma}}
		\quad\quad\forall t\geq 0;
		\label{th:decay-au-nc}
	\end{equation}
	\begin{equation}
		|\uep'(t)|^{2}\leq
		\frac{C_{2}}{(1+t)^{[2\gamma^{2}+(1-p)\gamma+p+1]/
		(\gamma^{2}+\gamma)}}
		\quad\quad\forall t\geq 0.
		\label{th:decay-u'-nc}
	\end{equation}
\end{thm}

The last result of this paper concerns the case $p>1$.  An analogous
result holds true for nondegenerate equations (see
\cite[Theorem~2.3]{gg:cwd}).

\begin{thm}\label{thm:p>1}
	Let $H$ and $A$ be as in Theorem~\ref{thm:nc-existence}. 
	Let $m:[0,+\infty)\to[0,+\infty)$ be a continuous function. Let
	$b:[0,+\infty)\to(0,+\infty)$ be a continuous function such that
	\begin{equation}
		\int_{0}^{+\infty}b(s)\,ds<+\infty.
		\label{hp:int-b}
	\end{equation}
	
	Let $(u_{0},u_{1})\in\da\times\dau$ be such that
	\begin{equation}
		|u_{1}|^{2}+\int_{0}^{|A^{1/2}u_{0}|^{2}}m(\sigma)\,d\sigma>0.
		\label{hp:p>1}
	\end{equation}
	
	Let us assume that for some $\ep>0$ problem (\ref{pbm:h-eq-gen}),
	(\ref{pbm:h-data}) has a global solution $\uep$ satisfying
	(\ref{th:uep-reg}).
	
	Then
	\begin{equation}
		\liminf_{t\to +\infty}\left(
		|\uep'(t)|^{2}+\auq{\uep(t)}\right)>0.
		\label{th:p>1}
	\end{equation}
\end{thm}

\begin{rmk}
	\begin{em}
		The constants $\ep_{0}$, $C_{1}$, $C_{2}$ given in
		Theorem~\ref{thm:c-existence} above may be taken as continuous
		functions of $\nu$, $\gamma$, $p$, $|u_{0}|$,
		$|A^{1/2}u_{0}|$, $|A^{1/2}u_{0}|^{-1}$, $|Au_{0}|$,
		$|u_{1}|$, $|A^{1/2}u_{1}|$. The same is true for the 
		constants $\ep_{0}$, $C_{1}$, $C_{2}$ given in
		Theorem~\ref{thm:nc-existence}, apart from the fact that in
		this case there is no dependence on $\nu$.
	\end{em}
\end{rmk}

\begin{rmk}
	\begin{em}
		Our results can be easily extended to more general Kirchhoff
		equations.  For example with the same technique we can deal
		with nonlinearities $m:[0,+\infty)\to[0,+\infty)$ of class
		$C^{1}$ such that 
		$$c_{1}\sigma^{\gamma}\leq m(\sigma)\leq c_{2}\sigma^{\gamma},
		\hspace{3em} c_{1}\sigma^{\gamma-1}\leq m'(\sigma)\leq
		c_{2}\sigma^{\gamma-1}$$
		in a right-hand neighborhood of $\sigma=0$ for suitable
		positive constants $c_{1}$ and $c_{2}$.  However this
		generality only complicates proofs without introducing any new
		idea.  
	\end{em}
\end{rmk}

\begin{rmk}
	\begin{em}
		In Theorem~\ref{thm:c-existence} we assume that $\gamma>0$,
		while in Theorem~\ref{thm:nc-existence} we assume that
		$\gamma\geq 1$.  Some weaker results can be obtained with
		similar techniques also when the operator is noncoercive and
		$\gamma>0$.  For example for every $\gamma>0$ one can prove
		the global solvability for every $p\in[0,\gamma/(\gamma+2)]$
		(and of course $\ep$ small enough).  The solution also
		satisfies (\ref{th:decay-a1/2u-nc}).  We sketch the argument
		in Remark~\ref{rmk:proof}.  Note that when $\gamma\geq 1$ the
		upper bound $\gamma/(\gamma+2)$ is always less than the upper
		bound in (\ref{hp:pnc}).
	\end{em}
\end{rmk}

\setcounter{equation}{0}
\section{Proofs}\label{sec:proofs}

Proofs are organized as follows.  First of all in~\ref{sec:ODE} we
state and prove two simple comparison results for ordinary
differential equations, which we need several times in the sequel.
Then we prove Theorem~\ref{thm:c-existence} and
Theorem~\ref{thm:nc-existence}.  Their proofs have a common part,
which we concentrate in~\ref{sec:common} in the form of an a priori
estimate (Proposition~\ref{prop:apriori}).  Then in~\ref{sec:c} we
conclude the proof of Theorem~\ref{thm:c-existence}, and
in~\ref{sec:nc} we conclude the proof of
Theorem~\ref{thm:nc-existence}.  Finally, in~\ref{sec:p>1} we prove
Theorem~\ref{thm:p>1}.

\subsection{Comparison results for ODEs}\label{sec:ODE}

Numerous variants of the following comparison result have already been
used in
\cite{ghisi1,g:non-lip,gg:k-dissipative,gg:k-decay,gg:k-PS,gg:cwd}.

\begin{lemma}\label{lemma:ode}
	Let $T>0$, let $p\geq 0$, and let $f:[0,T]\to[0,+\infty)$ be a
	function of class $C^{1}$.  Let us assume that there exist two
	constants $c_{1}>0$, $c_{2}\geq 0$  such that
	\begin{equation}
		f'(t)\leq -\frac{c_{1}}{(1+t)^{p}}f(t)+c_{2}\sqrt{f(t)}
		\hspace{2em}
		\forall t\in [0,T].
		\label{hp:lemma-ode}
	\end{equation}
	
	Then we have that
	\begin{equation}
		f(t)\leq f(0)+\left(\frac{c_{2}}{c_{1}}\right)^{2}(1+t)^{2p}
		\hspace{2em}
		\forall t\in [0,T].
		\label{th:lemma-ode}
	\end{equation} 
\end{lemma}

\prf From (\ref{hp:lemma-ode}) if follows that
$$f'(t)\leq -\frac{c_{1}}{2(1+t)^{p}}f(t)+
\frac{c_{2}^{2}}{2c_{1}}(1+t)^{p},$$
which is equivalent to say that $f(t)$ is a subsolution of the
differential equation
\begin{equation}
	y'(t)= -\frac{c_{1}}{2(1+t)^{p}}y(t)+
	\frac{c_{2}^{2}}{2c_{1}}(1+t)^{p}.
	\label{eq:ode}
\end{equation}

Let $g(t)$ denote the right-hand side of (\ref{th:lemma-ode}). Then
$$-\frac{c_{1}}{2(1+t)^{p}}g(t)+
\frac{c_{2}^{2}}{2c_{1}}(1+t)^{p} \leq 0\leq g'(t),$$
which is equivalent to say that $g(t)$ is a supersolution of
(\ref{eq:ode}).  Since $f(0)\leq g(0)$ the conclusion follows from the
standard comparison principle between subsolutions and supersolutions.
\qed

The next comparison result looks quite technical, but the basic idea
is the following.  When $f(t)\equiv 0$ the differential inequalities
(\ref{hp:ode-wy1}) and (\ref{hp:ode-wy2}) can be explicitly
integrated, thus providing estimates for $w(t)$.  When $f(t)$ is small
enough according to (\ref{hp:c1c2}) estimates of the same order can
still be proved.  The interested reader is referred to
\cite[Lemma~4.2]{gg:k-decay} for a similar comparison result.

\begin{lemma}\label{lemma:ode-vs}
	Let $p\geq 0$, $\gamma>0$, $\alpha>0$, and $T>0$ be real
	numbers.  Let $w:[0,T]\to[0,+\infty)$ be a function of class
	$C^{1}$ with $w(0)>0$, and let $f:[0,T]\to\re$ be a continuous
	function.
	
	Let us assume that 
	\begin{equation}
		\left|\int_{0}^{t}(1+s)^{p}f(s)\,ds\right|\leq
		\min\left\{\frac{1}{4\gamma[w(0)]^{\gamma}},
		\frac{\alpha}{2(p+1)}\right\}(1+t)^{p+1} \quad\quad\forall
		t\in[0,T].
		\label{hp:c1c2}
	\end{equation}
	
	Then we have the following implications.
	\begin{enumerate}
		\renewcommand{\labelenumi}{(\arabic{enumi})}
		\item If $w$ satisfies the differential inequality
		\begin{equation}
			w'(t)\leq -2(1+t)^{p}\left[w(t)\right]^{1+\gamma}
			\left(\alpha+f(t)\right)
			\quad\quad\forall t\in[0,T],
			\label{hp:ode-wy1}
		\end{equation}
		then we have the following estimate
		\begin{equation}
			w(t)\leq w(0)\left[\max\left\{
			2,\frac{p+1}{\alpha\gamma[w(0)]^{\gamma}}\right\}
			\right]^{1/\gamma}\cdot\frac{1}{(1+t)^{(p+1)/\gamma}}
			\quad\quad\forall t\in[0,T].
			\label{th:ode-wy1}
		\end{equation}
	
		\item If $w$ satisfies the differential inequality 
		\begin{equation}
			w'(t)\geq -2(1+t)^{p}\left[w(t)\right]^{1+\gamma}
			\left(\alpha+f(t)\right)
			\quad\quad\forall t\in[0,T],
			\label{hp:ode-wy2}
		\end{equation}
		then we have the following estimate
		\begin{equation}
			w(t)\geq w(0)\left[1+\frac{3\alpha\gamma
			\left[w(0)\right]^{\gamma}}{p+1}
			\right]^{-1/\gamma}\cdot\frac{1}{(1+t)^{(p+1)/\gamma}}
			\quad\quad\forall t\in[0,T].
			\label{th:ode-wy2}
		\end{equation}
	\end{enumerate}
\end{lemma}

\prf 
Let $y(t)$ be the solution of the Cauchy problem
$$y'(t)=-2\left[y(t)\right]^{\gamma+1},
\quad\quad
y(0)=w(0).$$

It is easy to see that
$$y(t)=w(0)\left(\strut 1+2\gamma[w(0)]^{\gamma}t\right)^{-1/\gamma}
\quad\quad
\forall t>-\frac{1}{2\gamma[w(0)]^{\gamma}}.$$

For every $t\in[0,T]$ let us set 
$$\Phi(t):=\frac{\alpha}{p+1}\left[(1+t)^{p+1}-1\right]+
\int_{0}^{t}(1+s)^{p}f(s)\,ds,
\hspace{3em}
z(t):=y\left(\Phi(t)\right).$$

First of all we have to prove that $z(t)$ is well defined.  This is
true because if we set 
$$C:=\min\left\{\frac{1}{4\gamma[w(0)]^{\gamma}},
\frac{\alpha}{2(p+1)}\right\},$$
then from (\ref{hp:c1c2}) we have that
$$\Phi(t)\geq
\left(\frac{\alpha}{p+1}-C\right)\left[(1+t)^{p+1}-1\right]-C\geq
-C>-\frac{1}{2\gamma[w(0)]^{\gamma}}.$$

Moreover a simple calculation shows that $z(t)$ is a solution of the
Cauchy problem
\begin{equation}
	z'(t)=-2(1+t)^{p}\left[z(t)\right]^{\gamma+1}(\alpha+f(t))
	\quad\quad\forall t\in[0,T],
	\label{ode:cp-eq}
\end{equation}
\begin{equation}
	z(0)=w(0).
	\label{ode:cp-data}
\end{equation}

\subparagraph{\textmd{\emph{Proof of statement (1)}}}

Assumption (\ref{hp:ode-wy1}) is equivalent to say that $w(t)$ is a
subsolution of the Cauchy problem (\ref{ode:cp-eq}),
(\ref{ode:cp-data}). The usual comparison principle implies that
$w(t)\leq z(t)$. Now we have to estimate $z(t)$. From (\ref{hp:c1c2}) 
it follows that
\begin{eqnarray*}
	1+2\gamma[w(0)]^{\gamma}\Phi(t) & \geq &
	1+2\gamma[w(0)]^{\gamma}\left(
	\left(\frac{\alpha}{p+1}-C\right)\left[(1+t)^{p+1}-1\right]
	-C\right) \\
	 & \geq & 1+2\gamma[w(0)]^{\gamma}\left(
	\frac{\alpha}{2(p+1)}\left[(1+t)^{p+1}-1\right]
	-\frac{1}{4\gamma[w(0)]^{\gamma}}\right)  \\
	 & = & \frac{1}{2}+ \frac{\alpha\gamma[w(0)]^{\gamma}}{p+1}
	 \left[(1+t)^{p+1}-1\right] \\
	 & \geq & \min\left\{\frac{1}{2},
	  \frac{\alpha\gamma[w(0)]^{\gamma}}{p+1}\right\}(1+t)^{p+1},
\end{eqnarray*}
where in the last step we exploited the elementary inequality
$$A+B(x-1)\geq\min\{A,B\}x
\quad\quad
\forall A\geq 0,\ \forall B\geq 0,\ \forall x\geq 1.$$

It follows that
\begin{eqnarray*}
	w(t)\;\leq\;z(t) & = &  w(0)\left[
	1+2\gamma[w(0)]^{\gamma}\Phi(t)\right]^{-1/\gamma}   \\
	 & \leq & w(0)\left[\max\left\{
	 2,\frac{p+1}{\alpha\gamma[w(0)]^{\gamma}}\right\}
	 \right]^{1/\gamma}\cdot\frac{1}{(1+t)^{(p+1)/\gamma}},
\end{eqnarray*}
which is exactly (\ref{th:ode-wy1}).

\subparagraph{\textmd{\emph{Proof of statement (2)}}}

Assumption (\ref{hp:ode-wy2}) is equivalent to say that $w(t)$ is a
supersolution of the Cauchy problem (\ref{ode:cp-eq}),
(\ref{ode:cp-data}), hence $w(t)\geq z(t)$ for every $t\in[0,T]$.

Since
\begin{eqnarray*}
	1+2\gamma[w(0)]^{\gamma}\Phi(t) & \leq &
	(1+t)^{p+1}+2\gamma[w(0)]^{\gamma}\left(
	\frac{\alpha}{p+1}+C\right)(1+t)^{p+1} \\
	 & \leq & \left(1+\frac{3\alpha\gamma[w(0)]^{\gamma}}{p+1}\right)
	 (1+t)^{p+1},
\end{eqnarray*}
the conclusion follows as in the previous case.
\qed

\subsection{Basic energy estimates}\label{sec:common}

In this section we prove some energy estimates and a lower bound for
$|A^{1/2}\uep(t)|$.  Such estimates don't require the coerciveness of
the operator, and they are fundamental both in the proof of
Theorem~\ref{thm:c-existence} and in the proof of
Theorem~\ref{thm:nc-existence}.  They extend to the weakly dissipative
equation the estimates stated in \cite[section~3.4]{gg:k-decay} in the
case of constant dissipation.

The estimates involve the following energies
\begin{equation}
	\Fep(t) :=  \ep\frac{|A^{1/2}\uep'(t)|^{2}}{\au{\uep(t)}^{2\gamma}}+
	|A\uep(t)|^{2},
	\label{defn:E} 
\end{equation}
\begin{equation}
	\Pep(t) :=  \ep
	\frac{\au{\uep(t)}^{2}\au{\uep'(t)}^{2}-
	\langle A\uep(t),\uep'(t)\rangle^{2}}{\au{\uep(t)}^{2\gamma+4}}+
	\frac{|A\uep(t)|^{2}}{\au{\uep(t)}^{2}},
	\label{defn:P} 
\end{equation}
\begin{equation}
	\Qep(t)  :=  \frac{|\uep'(t)|^{2}}{\au{\uep(t)}^{4\gamma+2}},
	\label{defn:Q} 
\end{equation}
\begin{equation}
	\Rep(t)  :=  \ep\frac{\au{\uep'(t)}^{2}}{\au{\uep(t)}^{2\gamma+2}}+
	\frac{|A\uep(t)|^{2}}{\au{\uep(t)}^{2}}.
	\label{defn:R} 
\end{equation}

We point out that the first summand in the definition of $\Pep(t)$ is 
nonnegative due to the Cauchy-Schwarz inequality.
	
We state the result in the form of an a priori estimate.  We assume
that in some interval $[0,S)$ there exists a solution of the hyperbolic
problem satisfying a given estimate (see (\ref{hp:K}) below), and we
deduce that this solution satisfies several energy inequalities in the
same interval.  We point out that all constants do not depend on $S$.

\begin{prop}\label{prop:apriori}
	Let $H$ and $A$ be as in Theorem~\ref{thm:nc-existence}. Let
	$\gamma>0$, $p\in[0,1]$, $K>0$ be real numbers, and let
	$(u_{0},u_{1})\in\da\times\dau$ satisfy (\ref{hp:mdg}). 
	
	Then there exists positive constants $\ep_{0}$, $\sigma_{0}$,
	$\sigma_{1}$ with the following property.  If $\ep\in(0,\ep_{0})$,
	$S>0$, and $$\uep\in C^{2}([0,S);H)\cap C^{1}([0,S);\dau) \cap
	C^{0}([0,S);\da)$$
	is a solution of (\ref{pbm:h-eq}), (\ref{pbm:h-data}) such that
	\begin{equation}
		A^{1/2}\uep(t)\neq 0 \mbox{ and }\frac{|\langle
		A\uep(t),\uep'(t)\rangle|}{\auq{\uep(t)}}
		\leq\frac{K}{(1+t)^{p}} \quad\quad\forall t\in[0,S),
		\label{hp:K}
	\end{equation}
	then for every $t\in[0,S)$ we have that
	\begin{equation}
		\Fep(t)+\int_{0}^{t}\frac{1}{(1+s)^{p}}
		\frac{\auq{\uep'(s)}}{\au{\uep(s)}^{2\gamma}}\,ds\leq
		\Fep(0);
		\label{th:est-F}
	\end{equation}
	\begin{equation}
		\Pep(t)\leq\Pep(0);
		\label{th:est-P}
	\end{equation}
	\begin{equation}
		\Qep(t)\leq\Qep(0)+4\Pep(0)(1+t)^{2p};
		\label{th:est-Q}
	\end{equation}
	\begin{eqnarray}
		\lefteqn{\hspace{-6em}
		(1+t)^{2p}\Rep(t)+\int_{0}^{t}(1+s)^{p}
		\frac{\auq{\uep'(s)}}{\au{\uep(s)}^{2\gamma+2}}\,ds
		\leq}
		\nonumber  \\
		\mbox{\hspace{2em}} & \leq &
		\left[\Rep(0)+2(K+1)\Pep(0)\right](1+t)^{p+1};
		\label{th:est-R}
	\end{eqnarray}
	\begin{equation}
		\left|\int_{0}^{t}(1+s)^{p}\frac{\langle \uep''(s),A\uep(s)\rangle}
		{\au{\uep(s)}^{2\gamma+2}}\,ds\right|\leq
		\sigma_{0}(1+t)^{p+1};
		\label{th:est-G}
	\end{equation}
	\begin{equation}
		\auq{\uep(t)}\geq
		\frac{\sigma_{1}}{(1+t)^{(p+1)/\gamma}}.
		\label{th:est-la1/2}
	\end{equation}
\end{prop}

\prf
Let us set 
$$ \sigma_{0} := \frac{|\langle u_{1},Au_{0}\rangle|}
{\au{u_{0}}^{2\gamma+2}}+\frac{3}{2}
\left(\sqrt{P_{1}(0)Q_{1}(0)}+2P_{1}(0)\right)+(2\gamma+3)\left(R_{1}(0)+
2(K+1)P_{1}(0)\right),$$
\begin{equation}
	\sigma_{1}:=\auq{u_{0}}\left(1+
	\frac{3\gamma P_{1}(0)|A^{1/2}u_{0}|^{2\gamma}}
	{p+1}\right)^{-1/\gamma}.
	\label{defn:sigma1}
\end{equation}

Let us choose $\ep_{0}$ in such a way that
\begin{equation}
	4\ep_{0}\leq 1,
	\quad\quad
	4\ep_{0}K(\gamma+1)\leq 1,
	\quad\quad
	\sigma_{0}\ep_{0}\leq\min\left\{
	\frac{1}{4\gamma\au{u_{0}}^{2\gamma}},
	\frac{P_{1}(0)}{2(p+1)} \right\}.
	\label{defn:ep0}
\end{equation}

\subparagraph{\textmd{\emph{Proof of (\ref{th:est-F}) through
(\ref{th:est-R})}}}

Let us compute the time derivative of the energies (\ref{defn:E})
through (\ref{defn:R}).  After some computations we find that
\begin{equation}
	\Fep'  = 
	-2\left(\frac{1}{(1+t)^{p}}+
	\gamma\ep\frac{\langle\uep',A\uep\rangle}{\au{\uep}^{2}}\right)
	\frac{|A^{1/2}\uep'|^{2}}{\au{\uep}^{2\gamma}},  
	\label{deriv:F}  
\end{equation}
\begin{equation}
	\Pep'  =  -2\left(\frac{1}{(1+t)^{p}}+
	(\gamma+2)\ep\frac{\langle\uep',A\uep\rangle}{\au{\uep}^{2}}\right)
	\frac{\au{\uep}^{2}\au{\uep'}^{2}-
	\langle A\uep,\uep'\rangle^{2}}{\au{\uep}^{2\gamma+4}},
	\label{deriv:P}  
\end{equation}
\begin{equation}
	\Qep'  =  -\frac{2}{\ep}\left(\frac{1}{(1+t)^{p}}+
	(2\gamma+1)
	\ep\frac{\langle\uep',A\uep\rangle}{\au{\uep}^{2}}\right)\Qep-
	\frac{2}{\ep}\frac{\langle\uep',A\uep\rangle}{\au{\uep}^{2\gamma+2}},
	\label{deriv:Q}  
\end{equation}
\begin{equation}
	\Rep'  =  -2\left(\frac{1}{(1+t)^{p}}+
	(\gamma+1)\ep\frac{\langle\uep',A\uep\rangle}{\au{\uep}^{2}}\right)
	\frac{\au{\uep'}^{2}}{\au{\uep}^{2\gamma+2}}-
	2\frac{\langle\uep',A\uep\rangle|A\uep|^{2}}{\au{\uep}^{4}}.
	\label{deriv:R}  
\end{equation}

Thanks to assumption (\ref{hp:K}) and the second inequality in
(\ref{defn:ep0}) we have that
\begin{equation}
	\Fep'(t)\leq -\frac{1}{(1+t)^{p}}
	\frac{|A^{1/2}\uep'(t)|^{2}}{\au{\uep(t)}^{2\gamma}};
	\label{est:DF}
\end{equation}
\begin{equation}
	\Pep'(t)\leq 0;
	\label{est:DP}
\end{equation}
\begin{equation}
	\Qep'(t)\leq -\frac{1}{\ep}\frac{1}{(1+t)^{p}}\Qep(t)-
	\frac{2}{\ep}\frac{\langle\uep'(t),A\uep(t)\rangle}
	{\au{\uep(t)}^{2\gamma+2}};
	\label{est:DQ}
\end{equation}
\begin{equation}
	\Rep'(t)\leq -\frac{3}{2}\frac{1}{(1+t)^{p}}
	\frac{|A^{1/2}\uep'(t)|^{2}}{\au{\uep(t)}^{2\gamma+2}}-
	2\frac{\langle\uep'(t),A\uep(t)\rangle|A\uep(t)|^{2}}
	{\au{\uep(t)}^{4}}.
	\label{est:DR}
\end{equation}

Integrating (\ref{est:DF}) in $[0,t]$ we obtain (\ref{th:est-F}).

Conclusion (\ref{th:est-P}) trivially follows from (\ref{est:DP}).

From (\ref{est:DQ}) we deduce that
\begin{eqnarray*}
	\Qep'(t) & \leq & -\frac{1}{\ep}\frac{1}{(1+t)^{p}}\Qep(t)+
	\frac{2}{\ep}\frac{|\uep'(t)|\cdot|A\uep(t)|}
	{\au{\uep(t)}^{2\gamma+2}}  \\
	 & \leq & -\frac{1}{\ep}\frac{1}{(1+t)^{p}}\Qep(t)+
	\frac{2}{\ep}\sqrt{\Pep(0)}\sqrt{\Qep(t)}.
\end{eqnarray*}

Therefore applying Lemma~\ref{lemma:ode} we obtain (\ref{th:est-Q}).

From (\ref{est:DR}) we have that
\begin{eqnarray*}
	\left[(1+t)^{2p}\Rep(t)\right]' & = & 
	2p(1+t)^{2p-1}\Rep(t)+(1+t)^{2p}\Rep'(t)   \\
	 & \leq  & 2p(1+t)^{2p-1}\ep
	 \frac{\au{\uep'}^{2}}{\au{\uep}^{2\gamma+2}}+
	2p(1+t)^{2p-1}\frac{|A\uep|^{2}}{\au{\uep}^{2}}+\\
	&&
	-\frac{3}{2}(1+t)^{p}
	\frac{|A^{1/2}\uep'|^{2}}{\au{\uep}^{2\gamma+2}}-
	2(1+t)^{2p}\frac{|A\uep|^{2}}{\au{\uep}^{2}}
	\frac{\langle\uep',A\uep\rangle}
	{\au{\uep}^{2}}  \\
	 & =: & I_{1}(t)+I_{2}(t)+I_{3}(t)+I_{4}(t).
\end{eqnarray*}

Since $2p-1\leq p$, and $2p\ep \leq 2\ep_{0}\leq 1/2$, we have that
$$I_{1}(t)+I_{3}(t)\leq \left(2p\ep-\frac{3}{2}\right)(1+t)^{p}
\frac{|A^{1/2}\uep'(t)|^{2}}{\au{\uep(t)}^{2\gamma+2}}\leq -(1+t)^{p}
\frac{|A^{1/2}\uep'(t)|^{2}}{\au{\uep(t)}^{2\gamma+2}}.$$

From (\ref{hp:K}), (\ref{th:est-P}), and the fact that $2p-1\leq p$ we
have that 
$$I_{2}(t)+I_{4}(t)\leq 2(K+p)(1+t)^{p}
\frac{|A\uep(t)|^{2}}{\au{\uep(t)}^{2}}\leq
2(K+1)(1+t)^{p}\Pep(0).$$

It follows that
$$\left[(1+t)^{2p}\Rep(t)\right]'\leq -(1+t)^{p}
\frac{|A^{1/2}\uep'(t)|^{2}}{\au{\uep(t)}^{2\gamma+2}}+
2(K+1)\Pep(0)(1+t)^{p}.$$

Integrating in $[0,t]$ we obtain (\ref{th:est-R}).

\subparagraph{\textmd{\emph{Proof of (\ref{th:est-G})}}}

Let us consider the following identity
\begin{eqnarray}
	(1+t)^{p}\frac{\langle
	\uep'',A\uep\rangle}{|A^{1/2}\uep|^{2\gamma+2}} & = & 
	\left[(1+t)^{p}\frac{\langle \uep',A\uep\rangle}
	{|A^{1/2}\uep|^{2\gamma+2}}\right]'-
	(1+t)^{p}\frac{|A^{1/2}\uep'|^{2}}{|A^{1/2}\uep|^{2\gamma+2}}+
	\nonumber \\
	 & & +(2\gamma+2)(1+t)^{p}\frac{\langle \uep',A\uep\rangle^{2}}
	 {|A^{1/2}\uep|^{2\gamma+4}}- p(1+t)^{p-1}\frac{\langle
	 \uep',A\uep\rangle}{|A^{1/2}\uep|^{2\gamma+2}}  
	 \nonumber \\
	  & =: & J_{1}(t)+J_{2}(t)+J_{3}(t)+J_{4}(t).
	  \label{defn:iiii}
\end{eqnarray}

In order to estimate the integral of the left-hand side, we estimate
the integrals of the four terms in the right-hand side.  By
(\ref{th:est-P}) and (\ref{th:est-Q}) we have that
\begin{eqnarray}
    \frac{|\langle\uep'(t),A\uep(t)\rangle|}
    {|A^{1/2}\uep(t)|^{2\gamma+2}} & \leq &
    \frac{|A\uep(t)|}{|A^{1/2}\uep(t)|}\cdot
    \frac{|\uep'(t)|}{|A^{1/2}\uep(t)|^{2\gamma+1}} 
    \nonumber\\
     & \leq & \sqrt{\Pep(0)}\sqrt{\Qep(0)+4\Pep(0)(1+t)^{2p}} 
     \nonumber\\
     & \leq & \left(\sqrt{\Pep(0)\Qep(0)}+2\Pep(0)\right)(1+t)^{p},
     \label{est:PQ}
\end{eqnarray}
hence
$$\left|\int_{0}^{t}J_{1}(s)\,ds\right|\leq
\frac{|\langle u_{1},Au_{0}\rangle|}
{|A^{1/2}u_{0}|^{2\gamma+2}}+(1+t)^{2p}
\left(\sqrt{\Pep(0)\Qep(0)}+2\Pep(0)\right).$$

The integral of $J_{2}(t)$ can be easily estimated using
(\ref{th:est-R}).

As for $J_{3}(t)$, by Cauchy-Schwarz inequality we have that
$$\frac{\langle\uep'(t),A\uep(t)\rangle^{2}}
{|A^{1/2}\uep(t)|^{2\gamma+4}}\leq
\frac{|A^{1/2}\uep'(t)|^{2}|A^{1/2}\uep(t)|^{2}}
{|A^{1/2}\uep(t)|^{2\gamma+4}}
=\frac{|A^{1/2}\uep'(t)|^{2}}{|A^{1/2}\uep(t)|^{2\gamma+2}},$$
and therefore we reduce once again to (\ref{th:est-R}).

Finally, from (\ref{est:PQ}) we obtain that
$$\left|\int_{0}^{t}J_{4}(s)\,ds\right|\leq
\frac{1}{2} \left(\sqrt{\Pep(0)\Qep(0)}+2\Pep(0)\right)
(1+t)^{2p}.$$

Plugging all these estimates in (\ref{defn:iiii}), and recalling once
again that $1\leq(1+t)^{2p}\leq(1+t)^{p+1}$ for every $t\geq 0$, we
obtain (\ref{th:est-G}).

\subparagraph{\textmd{\emph{Proof of (\ref{th:est-la1/2})}}}

Let us set $\wep(t):=\auq{\uep(t)}$.  Then
\begin{equation}
	\wep'(t)=-2(1+t)^{p}\left[\wep(t)\right]^{\gamma+1}
	\left(\frac{|A\uep(t)|^{2}}{|A^{1/2}\uep(t)|^{2}}+ \ep\frac{\langle
	\uep''(t),A\uep(t)\rangle}{|A^{1/2}\uep(t)|^{2\gamma+2}} \right),
	\label{deriv:wep}
\end{equation}
hence by (\ref{th:est-P})
$$\wep'(t)\geq -2(1+t)^{p}\left[\wep(t)\right]^{\gamma+1}
\left(P_{1}(0)+ \ep\frac{\langle
\uep''(t),A\uep(t)\rangle}{|A^{1/2}\uep(t)|^{2\gamma+2}} \right).$$

This means that $\wep$ satisfies a differential inequality of the form
(\ref{hp:ode-wy2}) with
\begin{equation}
	\alpha:=P_{1}(0),
	\hspace{3em}
	f(t):=\ep\frac{\langle
	\uep''(s),A\uep(s)\rangle}{|A^{1/2}\uep(s)|^{2\gamma+2}}.
	\label{defn:alpha-f}
\end{equation}

Thanks to (\ref{th:est-G}) and the last inequality in (\ref{defn:ep0})
we have that 
$$\left|\int_{0}^{t}(1+s)^{p}f(s)\,ds\right| \leq
\ep\sigma_{0}(1+t)^{p+1} \leq
\min\left\{\frac{1}{4\gamma|A^{1/2}u_{0}|^{2\gamma}},
\frac{P_{1}(0)}{2(p+1)}\right\}(1+t)^{p+1},$$
and therefore the function $f(t)$ satisfies assumption (\ref{hp:c1c2})
of Lemma~\ref{lemma:ode-vs}.  From statement~(2) of that lemma we
obtain (\ref{th:est-la1/2}).
%
\qed

\subsection{Proof in the coercive case}\label{sec:c}

\subparagraph{\textmd{\emph{Local maximal solutions}}}

Problem (\ref{pbm:h-eq}), (\ref{pbm:h-data}) admits a unique local-in-time
solution, and this solution can be continued to a solution defined in 
a maximal interval $[0,T)$, where either $T=+\infty$, or
\begin{equation}
	\limsup_{t\to T^{-}}\left(|A^{1/2}\uep'(t)|^{2}+
	|A\uep(t)|^{2}\right)=+\infty,
	\label{limsup}
\end{equation}
or
\begin{equation}
	\liminf_{t\to T^{-}}|A^{1/2}\uep(t)|^{2}=0.
	\label{liminf}
\end{equation}

We omit the proof of these standard results.  The interested reader is
referred to \cite{gg:k-dissipative} (see also \cite{AG}).

\subparagraph{\textmd{\emph{Preliminaries and notations}}}

Let $\nu$ satisfy (\ref{hp:coercive}), and let
$$\sigma_{2}:=|A^{1/2}u_{0}|^{2}\left[\max\left\{
2,\frac{p+1}{\nu\gamma|A^{1/2}u_{0}|^{2\gamma}}\right\}
\right]^{1/\gamma}.$$

Let $K$ be such that
\begin{equation}
	K>\frac{|\langle Au_{0},u_{1}\rangle|}{|A^{1/2}u_{0}|^{2}},
	\hspace{3em}
	K>\left(\sqrt{P_{1}(0)Q_{1}(0)}+2P_{1}(0)\right)
	\sigma_{2}^{\gamma}.
	\label{defn:K-c}
\end{equation}

Starting with this value of $K$ let us define $\sigma_{0}$ and
$\sigma_{1}$ as in the proof of Proposition~\ref{prop:apriori}, and
let us choose $\ep_{0}$ satisfying (\ref{defn:ep0}), and the further
requirement 
\begin{equation}
	\sigma_{0}\ep_{0}\leq\min\left\{
	\frac{\nu}{2(p+1)},\frac{1}{4\gamma|A^{1/2}u_{0}|^{2\gamma}}
	\right\}.
	\label{defn:ep0-bis}
\end{equation}

Let us finally set 
$$S:=\sup\left\{\tau\in[0,T): A^{1/2}\uep(t)\neq 0 \mbox{ and
}\frac{|\langle A\uep(t),\uep'(t)\rangle|}{\auq{\uep(t)}}
\leq\frac{K}{(1+t)^{p}}\ \forall t\in[0,\tau]\right\}.$$

From the mild nondegeneracy assumption (\ref{hp:mdg}) and the first
inequality in (\ref{defn:K-c}) it is easy to see that $S>0$.  Moreover
in the interval $[0,S)$ all the conclusions of
Proposition~\ref{prop:apriori} hold true.

\subparagraph{\textmd{\emph{Estimate from above for $\au{\uep(t)}$}}}

Let us set $\wep(t):=\auq{\uep(t)}$ as in the proof of
Proposition~\ref{prop:apriori}.  Once again $\wep(t)$ is a solution of
(\ref{deriv:wep}).  Since we are in the coercive case we have that
$|A\uep(t)|^{2}\geq\nu\auq{\uep(t)}$.  Therefore from
(\ref{deriv:wep}) it follows that 
$$\wep'(t)\leq
-2(1+t)^{p}\left[\wep(t)\right]^{\gamma+1} \left(\nu+ \ep\frac{\langle
\uep''(t),A\uep(t)\rangle}{|A^{1/2}\uep(t)|^{2\gamma+2}} \right),$$
which means that $\wep$ satisfies an inequality of the form
(\ref{hp:ode-wy1}) with $\alpha:=\nu$, and $f(t)$ defined as in
(\ref{defn:alpha-f}).  Thanks to (\ref{th:est-G}) and
(\ref{defn:ep0-bis}) the function $f(t)$ satisfies assumption
(\ref{hp:c1c2}) of Lemma~\ref{lemma:ode-vs}.  From statement~(1) of
that lemma we obtain that
\begin{equation}
	\auq{\uep(t)}\leq\frac{\sigma_{2}}{(1+t)^{(p+1)/\gamma}}
	\quad\quad\forall t\in[0,S).
	\label{th:est-ua1/2}
\end{equation}

\subparagraph{\textmd{\emph{Global existence}}}

We prove that $S=T=+\infty$. Let us assume by contradiction that $S<T$.
By definition of $S$ this means that
\begin{equation}
    \mbox{either}\quad
    \auq{\uep(S)}=0
    \quad\mbox{or}\quad
    \frac{|\langle A\uep(S),\uep'(S)\rangle|}{\auq{\uep(S)}}
    =\frac{K}{(1+S)^{p}}.
    \label{S-alternative}
\end{equation}

By continuity all the estimates proved so far hold true also for
$t=S$.  In particular (\ref{th:est-la1/2}) rules out the first
possibility in (\ref{S-alternative}).

%

%

On the other hand from (\ref{est:PQ}), (\ref{th:est-ua1/2}), and the
second inequality in (\ref{defn:K-c}), we have that
\begin{eqnarray*}
	\frac{|\langle A\uep(S),\uep'(S)\rangle|}{\auq{\uep(S)}} & \leq & 
	\frac{|A\uep(S)|}{|A^{1/2}\uep(S)|}\cdot
	\frac{|\uep'(S)|}{|A^{1/2}\uep(S)|^{2\gamma+1}}\cdot
	|A^{1/2}\uep(S)|^{2\gamma}   \\
	 & \leq & \left(\sqrt{P_{1}(0)Q_{1}(0)}+2P_{1}(0)\right)
	 (1+S)^{p}\cdot\frac{\sigma_{2}^{\gamma}}{(1+S)^{p+1}}  \\
	  & < & \frac{K}{1+S} \ \leq\  \frac{K}{(1+S)^{p}},
\end{eqnarray*}
which rules out the second possibility in (\ref{S-alternative}).

It remains to prove that $T=+\infty$.  Let us assume by contradiction
that $T<+\infty$.  Then the quoted local existence result says that
either (\ref{limsup}) or (\ref{liminf}) holds true.  

On the other hand now we know that (\ref{th:est-la1/2}) is satisfied
for every $t\in[0,T)$, which rules (\ref{liminf}) out.  Moreover from
(\ref{th:est-ua1/2}) we have that $\auq{\uep(t)}$ is uniformly bounded
from above in $[0,T)$, hence by (\ref{th:est-F}) it follows that also
$|A^{1/2}\uep'(t)|$ and $|A\uep(t)|$ are uniformly bounded from above
in $[0,T)$.  This rules (\ref{limsup}) out.

\subparagraph{\textmd{\emph{Decay estimates}}}

Let us prove estimates (\ref{th:decay-a1/2u}), (\ref{th:decay-au}),
and (\ref{th:decay-u'}).  Now we know that the solution is global, and
that all the estimates proved so far hold true for every $t\geq 0$.

Therefore (\ref{th:decay-a1/2u}) follows from (\ref{th:est-la1/2}) and
(\ref{th:est-ua1/2}).  Moreover from (\ref{th:est-P}) and the
coerciveness assumption
(\ref{hp:coercive}) we have that
$$\nu\leq\frac{|A\uep(t)|^{2}}{\auq{\uep(t)}}\leq P_{1}(0)
\quad\quad\forall t\geq 0,$$
hence (\ref{th:decay-au}) follows from (\ref{th:decay-a1/2u}).
Finally, (\ref{th:decay-u'}) follows from (\ref{th:est-Q}) and
(\ref{th:est-ua1/2}).
\qed

\subsection{Proof in the noncoercive case}\label{sec:nc}

\subparagraph{\textmd{\emph{Local maximal solutions}}}

As in the coercive case there exists a unique local-in-time solution
which can be continued to a solution defined in 
a maximal interval $[0,T)$, where either $T=+\infty$, or
(\ref{limsup}) holds true, or (\ref{liminf}) holds true.

\subparagraph{\textmd{\emph{Preliminaries and notations}}}

Let $\sigma_{1}$ be the constant defined in (\ref{defn:sigma1}), 
let
$$\sigma_{3}:=16(\gamma+1)\left(|u_{1}|^{2}+
\au{u_{0}}^{2\gamma+2}+2|u_{0}|^{2}\right),$$
$$\sigma_{4}:=2\frac{|A^{1/2}u_{1}|^{2}}{|A^{1/2}u_{0}|^{2\gamma}}+
2|Au_{0}|^{2}+\frac{1}{2}
\frac{|\langle Au_{0},u_{1}\rangle|}{\au{u_{0}}^{2\gamma}}+
36\sigma_{1}^{1-\gamma},$$
and let $K$ be such that
\begin{equation}
	K>\frac{|\langle Au_{0},u_{1}\rangle|}{\auq{u_{0}}},
	\hspace{2em}
	K>\left[(1+\gamma)\sigma_{3}\right]^{(\gamma-1)/(\gamma+1)}
	\left(\frac{|u_{1}|}{\au{u_{0}}^{2\gamma}}
	\sqrt{\sigma_{4}}+4\sigma_{4}\right).
	\label{defn:K-nc}
\end{equation}

Starting with this value of $K$ let us define $\sigma_{0}$ as in the
proof of Proposition~\ref{prop:apriori}, and let us choose $\ep_{0}$
satisfying (\ref{defn:ep0}) and the further condition
$$16\ep_{0}\leq 1.$$

As in the coercive case let us finally set 
$$S:=\sup\left\{\tau\in[0,T): A^{1/2}\uep(t)\neq 0 \mbox{ and
}\frac{|\langle A\uep(t),\uep'(t)\rangle|}{\auq{\uep(t)}}
\leq\frac{K}{(1+t)^{p}}\ \forall t\in[0,\tau]\right\}.$$

From the mild nondegeneracy assumption (\ref{hp:mdg}), and the first
inequality in (\ref{defn:K-nc}), it is easy to see that $S>0$.
Moreover in the interval $[0,S)$ all the conclusions of
Proposition~\ref{prop:apriori} hold true.

In the following we set
$$\beta=\frac{p+1}{\gamma},$$
and we prove estimates involving the following energies
\begin{equation}
	\Hep(t):=\ep|\uep'(t)|^{2}+\frac{1}{\gamma+1}\au{\uep(t)}^{2\gamma+2};
	\label{defn:hep}
\end{equation}
\begin{equation}
	\Dep(t):=\ep(1+t)^{p}\langle\uep'(t),\uep(t)\rangle+
	\frac{1}{2}\left(1-\frac{\ep p}{(1+t)^{1-p}}\right)|\uep(t)|^{2};
	\label{defn:dep}
\end{equation}
\begin{equation}
	\DDep(t):=\ep(1+t)^{2\beta-1}\frac{\langle\uep'(t),A\uep(t)\rangle}
	{\au{\uep(t)}^{2\gamma}};
	\label{defn:ddep}
\end{equation}
\begin{equation}
	\Gep(t):=(1+t)^{\beta}\frac{|\uep'(t)|^{2}}{|A^{1/2}\uep(t)|^{4\gamma}}.
	\label{defn:gep}
\end{equation}

All the estimates we present are first claimed in the interval
$[0,S)$.  At the end of the proof we show that $S=T=+\infty$, thus
obtaining that all the estimates actually hold true for every $t\geq
0$.

\subparagraph{\textmd{\emph{First order estimate}}}

In this section of the proof we show that
\begin{equation}
	(1+t)^{p+1}\Hep(t)+|\uep(t)|^{2}+
	\int_{0}^{t}(1+s)|\uep'(s)|^{2}\,ds\leq\sigma_{3}
	\quad\quad\forall t\in[0,S).
	\label{est:hep}
\end{equation}

To this end we begin by taking the time derivative of
(\ref{defn:dep}):
$$\Dep'(t)=-(1+t)^{p}|A^{1/2}\uep(t)|^{2\gamma+2}+
\ep(1+t)^{p}|\uep'(t)|^{2}+
\frac{\ep p(1-p)}{2}\frac{|\uep(t)|^{2}}{(1+t)^{2-p}}.$$

Integrating in $[0,t]$ we obtain that
\begin{eqnarray}
	\int_{0}^{t}(1+s)^{p}|A^{1/2}\uep(s)|^{2\gamma+2}\,ds & = & 
	\Dep(0)-\Dep(t)+\ep\int_{0}^{t}(1+s)^{p}|\uep'(s)|^{2}\,ds+
	\nonumber  \\
	 &  & +\frac{\ep p(1-p)}{2}\int_{0}^{t}
	 \frac{|\uep(s)|^{2}}{(1+s)^{2-p}}\,ds.
	\label{est:0D}
\end{eqnarray}

From our assumptions on $\ep$ and $p$ we have that $2\ep<1/4$, $2p\leq
p+1$, $\ep p\leq 1/2$. Therefore
\begin{eqnarray*}
	-\Dep(t) & \leq & 2\ep^{2}(1+t)^{2p}|\uep'(t)|^{2}+
	\frac{1}{8}|\uep(t)|^{2}+\frac{\ep p}{2(1+t)^{1-p}}|\uep(t)|^{2}
	-\frac{1}{2}|\uep(t)|^{2}  \\ 
	& \leq &
	 \frac{1}{4}\ep(1+t)^{p+1}|\uep'(t)|^{2}-
	 \frac{1}{8}|\uep(t)|^{2}.
\end{eqnarray*}

Plugging this estimate in (\ref{est:0D}) we obtain that
\begin{eqnarray}
	\lefteqn{
	\hspace{-4em}
	\frac{1}{8}|\uep(t)|^{2}+
	\int_{0}^{t}(1+s)^{p}|A^{1/2}\uep(s)|^{2\gamma+2}\,ds\ \leq\ 
	\Dep(0)+\frac{1}{4}\ep(1+t)^{p+1}|\uep'(t)|^{2}+} 
	\nonumber  \\
	 & & +\ep\int_{0}^{t}(1+s)|\uep'(s)|^{2}\,ds+ \frac{\ep p(1-p)}{2}
	 \int_{0}^{t}\frac{|\uep(s)|^{2}}{(1+s)^{2-p}}\,ds.
	\label{est:D1}
\end{eqnarray}

Let us consider now the energy defined in (\ref{defn:hep}). A simple
calculation gives that
$$\left[(1+t)^{p+1}\Hep\right]'=-(1+t)\left(
2-\frac{\ep(p+1)}{(1+t)^{1-p}}\right)|\uep'|^{2}+
\frac{p+1}{\gamma+1}(1+t)^{p}|A^{1/2}\uep|^{2\gamma+2}.$$

Let us integrate in $[0,t]$.  Using (\ref{est:D1}) and rearranging the
terms we obtain that
\begin{eqnarray*}
	\lefteqn{
	\hspace{-4em} 
	(1+t)^{p+1}\left(1-\frac{1}{4}\frac{p+1}{\gamma+1}\right)
	\ep|\uep'(t)|^{2}+\frac{(1+t)^{p+1}}{\gamma+1}
	|A^{1/2}\uep(t)|^{2\gamma+2}
	\ \leq}  \\
	 \mbox{\hspace{3em}} & \leq & \Hep(0)-\left(
	 2-\ep(p+1)-\ep\frac{p+1}{\gamma+1}\right)
	 \int_{0}^{t}(1+s)|\uep'(s)|^{2}\,ds+ \\
	 & & +\frac{p+1}{\gamma+1}\left(\Dep(0)-\frac{1}{8}|\uep(t)|^{2}+
	 \frac{\ep p(1-p)}{2}
	 \int_{0}^{t}\frac{|\uep(s)|^{2}}{(1+s)^{2-p}}\,ds\right).
\end{eqnarray*}

From the smallness assumptions on $\ep$, and the fact that
$(p+1)/(\gamma+1)\leq 2$, it follows that
\begin{eqnarray}
	\lefteqn{ 
	\hspace{-4em} 
	\frac{1}{2}(1+t)^{p+1}\Hep(t)+
	\int_{0}^{t}(1+s)|\uep'(s)|^{2}\,ds+
	\frac{1}{8}\frac{p+1}{\gamma+1}|\uep(t)|^{2}\ \leq }
	\nonumber \\
	 & \leq &
	 \left(\Hep(0)+2\left|\Dep(0)\right|\right)+\frac{p+1}{\gamma+1}
	 \frac{\ep p(1-p)}{2}\int_{0}^{t}\frac{|\uep(s)|^{2}}{(1+s)^{2-p}}\,ds.
	\label{est:H1}
\end{eqnarray}

In particular we have that
$$|\uep(t)|^{2}\leq\frac{8(\gamma+1)}{p+1}
\left(\Hep(0)+2\left|\Dep(0)\right|\right)+4\ep(1-p)
\int_{0}^{t}\frac{|\uep(s)|^{2}}{(1+s)^{2-p}}\,ds,$$
hence by Gronwall's lemma 
\begin{eqnarray*}
	|\uep(t)|^{2} & \leq & \frac{8(\gamma+1)}{p+1}
	\left(\Hep(0)+2\left|\Dep(0)\right|\right)\exp\left(4\ep(1-p)
	\int_{0}^{t}\frac{1}{(1+s)^{2-p}}\,ds \right) \\
	 & \leq & \frac{8(\gamma+1)}{p+1}
	 \left(\Hep(0)+2\left|\Dep(0)\right|\right)\exp(4\ep)  \\
	  & \leq & 
	 \frac{16(\gamma+1)}{p+1}
	 \left(\Hep(0)+2\left|\Dep(0)\right|\right).
\end{eqnarray*}

Integrating in $[0,t]$ we obtain that
$$(1-p)\int_{0}^{t}\frac{|\uep(s)|^{2}}{(1+s)^{2-p}}\,ds\leq
16\frac{\gamma+1}{p+1} \left(\Hep(0)+2\left|\Dep(0)\right|\right).$$

Coming back to (\ref{est:H1}) we have therefore that
\begin{eqnarray}
	\lefteqn{ 
	\hspace{-8em} 
	\frac{1}{2}(1+t)^{p+1}\Hep(t)+
	\int_{0}^{t}(1+s)|\uep'(s)|^{2}\,ds+
	\frac{1}{8}\frac{p+1}{\gamma+1}|\uep(t)|^{2}}
	\nonumber  \\
	 & \leq & (1+8p\ep) \left(\Hep(0)+2\left|\Dep(0)\right|\right).
	\label{est:H2}
\end{eqnarray}

It remains to estimate the right-hand side. This can be easily done
because $8p\ep\leq 1$, and
\begin{eqnarray*}
	\Hep(0)+2\left|\Dep(0)\right| & \leq & 
	\ep|u_{1}|^{2}+\frac{1}{\gamma+1}|A^{1/2}u_{0}|^{2\gamma+2}+
	2\ep|\langle u_{1},u_{0}\rangle|+(1-\ep p)|u_{0}|^{2}\\
	 & \leq &
	 |u_{1}|^{2}+|A^{1/2}u_{0}|^{2\gamma+2}+
	 2|u_{0}|^{2}.
\end{eqnarray*}

Plugging this estimate in (\ref{est:H2}), and multiplying by
$8(\gamma+1)$, we obtain (\ref{est:hep}).

\subparagraph{\textmd{\emph{Second order estimate}}}

In this section of the proof we show that
\begin{equation}
	(1+t)^{\beta}\Fep(t)+\frac{1}{2}\frac{1}{(1+t)^{\beta}}
	\int_{0}^{t}(1+s)^{2\beta-p}
	\frac{\auq{\uep'(s)}}{\au{\uep(s)}^{2\gamma}}\,ds\leq\sigma_{4}
	\quad\quad\forall t\in[0,S).
	\label{est:fep}
\end{equation}

To this end we begin by computing the time derivative of
(\ref{defn:ddep}):
\begin{eqnarray}
	\DDep'(t) & = & -(1+t)^{2\beta-1}|A\uep(t)|^{2}+
	\ep(1+t)^{2\beta-1}\frac{|A^{1/2}\uep'(t)|^{2}}{|A\uep(t)|^{2\gamma}}+
	\nonumber  \\
	&&
	-2\gamma\ep(1+t)^{2\beta-1}\frac{\langle\uep'(t),A\uep(t)\rangle^{2}}
	{|A^{1/2}\uep(t)|^{2\gamma+2}}+
	\nonumber \\
	 &  & -(1+t)^{2\beta-p-1}\left(
	1-\ep\frac{2\beta-1}{(1+t)^{1-p}}\right)
	\frac{\langle\uep'(t),A\uep(t)\rangle}{\au{\uep(t)}^{2\gamma}}
	\nonumber \\
	&=:& I_{1}(t)+I_{2}(t)+I_{3}(t)+I_{4}(t).
	\label{deriv:ddep}
\end{eqnarray}

Let us estimate this derivative from above.  To this end in $I_{2}(t)$
we replace the exponent $2\beta-1$ with the bigger exponent
$2\beta-p$.  The term $I_{3}(t)$ is nonpositive and can be neglected.
In order to estimate $I_{4}(t)$ we remark that $0<\beta\leq 2$, hence
$0\leq|2\beta-1|\leq 3$.  Due to the smallness of $\ep$ we have
therefore that
$$\left|1-\ep\frac{2\beta-1}{(1+t)^{1-p}}\right|\leq
1+\frac{|2\beta-1|\ep}{(1+t)^{1-p}}\leq 1+3\ep\leq 2,$$
and thus
\begin{eqnarray*}
	\left|I_{4}(t)\right| & \leq & 2(1+t)^{2\beta-p-1}
	\frac{|A^{1/2}\uep'(t)|}{|A^{1/2}\uep(t)|^{2\gamma-1}} \\
	 & \leq & \frac{1}{4\beta}(1+t)^{2\beta-p}\frac{\auq{\uep'(t)}}
	{\au{\uep(t)}^{2\gamma}}+4\beta(1+t)^{2\beta-p-2}
	\frac{1}{\au{\uep(t)}^{2\gamma-2}}.
\end{eqnarray*}

Since $\gamma\geq 1$ we can estimate the last term using
(\ref{th:est-la1/2}).  After some calculations with the exponents we
obtain that
\begin{equation}
	\frac{1}{\au{\uep(t)}^{2\gamma-2}}\leq\sigma_{1}^{1-\gamma}
	(1+t)^{p+1-\beta},
	\label{est:inversa}
\end{equation}
hence
$$
\left|I_{4}(t)\right|
\leq\frac{1}{4\beta}(1+t)^{2\beta-p}\frac{\auq{\uep'(t)}}
{\au{\uep(t)}^{2\gamma}}+4\beta\sigma_{1}^{1-\gamma}
(1+t)^{\beta-1}.$$

Plugging these estimates in (\ref{deriv:ddep}) we have proved that
\begin{eqnarray*}
	\DDep'(t) & \leq & -(1+t)^{2\beta-1}|A\uep(t)|^{2}+
	(1+t)^{2\beta-p}\left(\ep+\frac{1}{4\beta}\right)\frac{\auq{\uep'(t)}}
	{\au{\uep(t)}^{2\gamma}}+ \\
	 &  & +4\beta\sigma_{1}^{1-\gamma}
(1+t)^{\beta-1}.
\end{eqnarray*}

Integrating in $[0,t]$ we obtain that
\begin{eqnarray}
	\int_{0}^{t}(1+s)^{2\beta-1}|A\uep(s)|^{2}\,ds & \leq &
	\left(\ep+\frac{1}{4\beta}\right)
	 \int_{0}^{t}(1+s)^{2\beta-p}\frac{\auq{\uep'(s)}}
	 {\au{\uep(s)}^{2\gamma}}\,ds+
	\nonumber  \\
	 & & +\DDep(0)-\DDep(t)+4\sigma_{1}^{1-\gamma} (1+t)^{\beta}.
	\label{est:DD1}
\end{eqnarray}

Using (\ref{est:inversa}) once more we have that
\begin{eqnarray*}
	-\DDep(t) & \leq & \frac{\ep^{2}}{2}(1+t)^{2\beta}
	\frac{\auq{\uep'(t)}}{\au{\uep(t)}^{2\gamma}}+
	\frac{1}{2}(1+t)^{2\beta-2}\frac{1}{\au{\uep(t)}^{2\gamma-2}}
	 \\
	 & \leq & \frac{\ep^{2}}{2}(1+t)^{2\beta}
	\frac{\auq{\uep'(t)}}{\au{\uep(t)}^{2\gamma}}+
	\frac{1}{2}\sigma_{1}^{1-\gamma}(1+t)^{\beta+p-1}.  
\end{eqnarray*}

Since $\beta+p-1\leq\beta$, plugging this estimate in (\ref{est:DD1})
we obtain that
\begin{eqnarray}
	\int_{0}^{t}(1+s)^{2\beta-1}|A\uep(s)|^{2}\,ds & \leq & 
	\DDep(0)+\frac{\ep^{2}}{2}(1+t)^{2\beta}
	\frac{\auq{\uep'(t)}}{\au{\uep(t)}^{2\gamma}}+
	\frac{9}{2}\sigma_{1}^{1-\gamma}(1+t)^{\beta}+
	\nonumber  \\
	 &  & +\left(\ep+\frac{1}{4\beta}\right)
	 \int_{0}^{t}(1+s)^{2\beta-p}\frac{\auq{\uep'(s)}}
	 {\au{\uep(s)}^{2\gamma}}\,ds.
	\label{est:DFH3}
\end{eqnarray}

Let us consider now the energy defined in (\ref{defn:E}). A simple
calculation gives that
\begin{eqnarray*}
	\left[(1+t)^{2\beta}\Fep\right]' & = & -(1+t)^{2\beta}\left(
	\frac{2}{(1+t)^{p}}+2\ep\gamma
	\frac{\langle A\uep,\uep'\rangle}{\auq{\uep}}-
	\frac{2\beta\ep}{1+t}\right)
	\frac{\auq{\uep'}}{\au{\uep}^{2\gamma}}+ \\
	 &  & +2\beta(1+t)^{2\beta-1}|A\uep|^{2}.
\end{eqnarray*}

By definition of $S$, and the second inequality in (\ref{defn:ep0}),
we have that
\begin{equation}
	2\ep\gamma
	\frac{\left|\langle A\uep(t),\uep'(t)\rangle\right|}{\auq{\uep(t)}}+
	\frac{2\beta\ep}{1+t}\leq
	\left(2\ep\gamma K+2\beta\ep\right)\frac{1}{(1+t)^{p}}\leq
	\frac{1}{(1+t)^{p}},
	\label{est:KS}
\end{equation}
hence
$$\left[(1+t)^{2\beta}\Fep(t)\right]'\leq
-(1+t)^{2\beta-p}
\frac{\auq{\uep'(t)}}{\au{\uep(t)}^{2\gamma}}
+2\beta(1+t)^{2\beta-1}|A\uep(t)|^{2}.$$

Let us integrate in $[0,t]$.  Using (\ref{est:DFH3}) and rearranging
the terms we obtain that
\begin{eqnarray*}
	(1+t)^{2\beta}\left(1-\beta\ep\right)\ep
	\frac{\auq{\uep'(t)}}{\au{\uep(t)}^{2\gamma}}+
	(1+t)^{2\beta}|A\uep(t)|^{2}+ &  &  \\
	+\left(\frac{1}{2}-2\beta\ep\right)\int_{0}^{t}(1+s)^{2\beta-p}
	\frac{\auq{\uep'(s)}}{\au{\uep(s)}^{2\gamma}}\,ds & \leq &
	\Fep(0)+2\beta\DDep(0)+  \\
	&& +9\beta\sigma_{1}^{1-\gamma}(1+t)^{\beta}.
\end{eqnarray*}

Since $\beta\leq 2$, and $2\beta\ep\leq 4\ep\leq 1/4$, it follows that
\begin{eqnarray*}
	\lefteqn{
	\hspace{-4em}
	\frac{1}{2}(1+t)^{2\beta}\Fep(t)+\frac{1}{4}\int_{0}^{t}
	(1+s)^{2\beta-p}
	\frac{\auq{\uep'(s)}}{\au{\uep(s)}^{2\gamma}}\,ds}  \\
	\mbox{\hspace{3em}} & \leq &
	 \Fep(0)+2\beta\DDep(0)+9\beta\sigma_{1}^{1-\gamma}(1+t)^{\beta}  
	 \\
	  & \leq & \frac{|A^{1/2}u_{1}|^{2}}{|A^{1/2}u_{0}|^{2\gamma}}+
	  |Au_{0}|^{2}+\frac{1}{4}\frac{|\langle
	  Au_{0},u_{1}\rangle|}{\au{u_{0}}^{2\gamma}}
	  +18\sigma_{1}^{1-\gamma}(1+t)^{\beta}.  
\end{eqnarray*}

Dividing by $(1+t)^{\beta}/2$ we obtain (\ref{est:fep}).

\subparagraph{\textmd{\emph{Estimate on the derivative}}}

Let us consider the energy defined in (\ref{defn:gep}). Its time
derivative is given by
\begin{eqnarray*}
	\Gep'(t) & = & -\frac{1}{\ep}(1+t)^{\beta}\left(
	\frac{2}{(1+t)^{p}}+4\ep\gamma
	\frac{\langle\uep'(t),A\uep(t)\rangle}{\auq{\uep(t)}}-
	\frac{\beta\ep}{1+t}\right)
	\frac{|\uep'(t)|^{2}}{\au{\uep(t)}^{4\gamma}}+ \\
	 &  & - \frac{2}{\ep}(1+t)^{\beta}
	\frac{\langle\uep'(t),A\uep(t)\rangle}{\au{\uep(t)}^{2\gamma}}.
\end{eqnarray*}

Arguing as in (\ref{est:KS}) we find that
$$4\ep\gamma
\frac{\left|\langle\uep'(t),A\uep(t)\rangle\right|}{\auq{\uep(t)}}+
\frac{\beta\ep}{1+t}\leq
\frac{3}{2}\frac{1}{(1+t)^{p}},$$
hence
$$\Gep'(t)\leq -\frac{1}{2\ep}\frac{1}{(1+t)^{p}}\Gep(t)+
\frac{2}{\ep}(1+t)^{\beta/2}|A\uep(t)|\cdot\sqrt{\Gep(t)}.$$

Thanks to (\ref{est:fep}) we have therefore that
$$\Gep'(t)\leq -\frac{1}{2\ep}\frac{1}{(1+t)^{p}}\Gep(t)+
\frac{2}{\ep}\sqrt{\sigma_{4}}\cdot\sqrt{\Gep(t)},$$
hence by Lemma~\ref{lemma:ode}
\begin{equation}
	\Gep(t)\leq\Gep(0)+16\sigma_{4}(1+t)^{2p}
	\quad\quad\forall t\in[0,S).
	\label{est:gep}
\end{equation}

\subparagraph{\textmd{\emph{Global existence}}}

We prove that $S=T=+\infty$.  Let us assume by contradiction that
$S<T$.  Then by continuity all the estimates
proved so far hold true also for $t=S$.  Moreover by definition of $S$
we have the alternative (\ref{S-alternative}).

The first possibility can be ruled out using (\ref{th:est-la1/2})
exactly as in the coercive case.

In order to rule out the second possibility we consider the inequality
\begin{equation}
	\frac{|\langle A\uep(S),\uep'(S)\rangle|}{\auq{\uep(S)}} \leq
	|A\uep(S)|\cdot\frac{|\uep'(S)|}{|A^{1/2}\uep(S)|^{2\gamma}}
	\cdot|A^{1/2}\uep(S)|^{2\gamma-2}.
	\label{est:prefinale}
\end{equation}

Let us estimate the three factors.  From (\ref{est:fep}) we have that
$$|A\uep(S)|\leq\frac{\sqrt{\sigma_{4}}}{(1+S)^{\beta/2}}.$$

From (\ref{est:gep}) we have that
$$\frac{|\uep'(S)|}{|A^{1/2}\uep(S)|^{2\gamma}}
\leq\frac{\sqrt{\Gep(0)+16\sigma_{4}(1+S)^{2p}}}
{(1+S)^{\beta/2}}\leq\left(\frac{|u_{1}|}{|A^{1/2}u_{0}|^{2\gamma}}+
4\sqrt{\sigma_{4}}\right)\frac{1}{(1+S)^{\beta/2-p}}.$$

Since $\gamma\geq 1$ the last factor in (\ref{est:prefinale}) can be
estimated using (\ref{est:hep}).  We obtain that
$$|A^{1/2}\uep(S)|^{2\gamma-2}\leq
\left[(\gamma+1)\sigma_{3}\right]^{(\gamma-1)/(\gamma+1)}
\frac{1}{(1+S)^{(p+1)(\gamma-1)/(\gamma+1)}}.$$

Plugging all these estimates in (\ref{est:prefinale}), and recalling
(\ref{defn:K-nc}), we obtain that
$$\frac{|\langle A\uep(S),\uep'(S)\rangle|}{\auq{\uep(S)}}<
\frac{K}{(1+S)^{p}}(1+S)^{2p-\beta-(p+1)(\gamma-1)/(\gamma+1)}.$$

If $p$ satisfies (\ref{hp:pnc}), then the last exponent is less than
or equal to zero, and this is enough to rule out the second
possibility in (\ref{S-alternative}).

It remains to prove that $T=+\infty$.  So let us assume by
contradiction that $T<+\infty$.  Then the quoted local existence
result says that either (\ref{limsup}) or (\ref{liminf}) holds true.

On the other hand as in the coercive case we have that
(\ref{th:est-la1/2}) is satisfied for every $t\in[0,T)$, which rules
(\ref{liminf}) out.  Moreover from (\ref{est:hep}) we have that
$\auq{\uep(t)}$ is uniformly bounded from above in $[0,T)$, hence by
(\ref{th:est-F}) it follows that also $|A^{1/2}\uep'(t)|$ and
$|A\uep(t)|$ are uniformly bounded from above in $[0,T)$.  This rules
(\ref{limsup}) out.

\subparagraph{\textmd{\emph{Decay estimates}}}

Let us prove estimates (\ref{th:decay-a1/2u-nc}), (\ref{th:decay-au-nc}),
and (\ref{th:decay-u'-nc}). Now we know that the solution is global,
and that all the estimates proved so far hold true for every $t\geq 0$.

Therefore the lower bound in (\ref{th:decay-a1/2u-nc}) follows from
(\ref{th:est-la1/2}), while the upper bound follows from
(\ref{est:hep}).  Moreover (\ref{th:decay-au-nc}) follows from
(\ref{est:fep}).  Finally, (\ref{th:decay-u'-nc}) follows from
(\ref{est:gep}) and (\ref{th:decay-a1/2u-nc}).
\qed

\begin{rmk}\label{rmk:proof}
	\begin{em}
		A careful inspection of the proofs reveals that
		(\ref{est:hep}) was proved without using the assumption
		$\gamma\geq 1$.  At this point one can modify
		(\ref{est:prefinale}) as follows: 
		$$\frac{|\langle A\uep(S),\uep'(S)\rangle|}{\auq{\uep(S)}}
		\leq
		\frac{|A\uep(S)|}{|A^{1/2}\uep(S)|}\cdot
		\frac{|\uep'(S)|}{|A^{1/2}\uep(S)|^{2\gamma+1}}
		\cdot|A^{1/2}\uep(S)|^{2\gamma}.$$
		
		Now we can estimate the first and second factor using
		(\ref{est:PQ}) as we did in the coercive case, and then
		estimate the last factor using (\ref{est:hep}).  All these
		inequalities require neither the coerciveness of the operator,
		nor $\gamma\geq 1$.
		
		We end up with an estimate such as
		$$\frac{|\langle A\uep(S),\uep'(S)\rangle|}{\auq{\uep(S)}}	
		\leq\frac{K_{1}}{(1+S)^{p}}(1+S)^{2p-\gamma(p+1)/
		(\gamma+1)}$$
		for a suitable constant $K_{1}$.  The last exponent is less
		than or equal to zero provided that $p\leq\gamma/(\gamma+2)$.
		This is the key point of the proof of global solvability for
		$\gamma>0$ and $p\in[0,\gamma/(\gamma+2)]$ without coerciveness
		assumptions.  We leave the details to the interested reader.
	\end{em}
\end{rmk}

\subsection{Proof of Theorem~\ref{thm:p>1}}\label{sec:p>1}

In analogy with (\ref{defn:hep}) let us set
$$\Hep(t):=\ep|\uep'(t)|^{2}+
\int_{0}^{|A^{1/2}\uep(t)|^{2}}m(\sigma)\,d\sigma.$$

Assumption (\ref{hp:p>1}) is equivalent to say that $\Hep(0)>0$.
Moreover we have that
$$\Hep'(t)=-2b(t)|\uep'(t)|^{2}\geq -\frac{2}{\ep}b(t)\Hep(t)
\quad\quad\forall t\geq 0,$$
hence
$$\Hep(t)\geq \Hep(0)\exp\left(-\frac{2}{\ep}\int_{0}^{t}b(s)\,ds\right)
\quad\quad\forall t\geq 0.$$

The right-hand side is greater than a positive constant independent on
$t$ because of (\ref{hp:int-b}) and the fact that $\Hep(0)>0$.  This
implies (\ref{th:p>1}).
\qed

%

\label{NumeroPagine}


\begin{thebibliography}{99}

\bibitem{AG}{\sc A.\ Arosio, S.\ Garavaldi}; On the mildly degenerate
Kirchhoff string.  \emph{Math.\  Methods Appl.\  Sci.}\ \textbf{14}
(1991), no.~3, 177--195.

\bibitem{debrito}{\sc E.\ H.\ de Brito}; The damped elastic stretched
string equation generalized: existence, uniqueness, regularity and
stability.  \emph{Applicable Anal.}\ \textbf{13} (1982), no.~3,
219--233.


\bibitem{ew}{\sc B.\ F.\ Esham, R.\ J.\ Weinacht}; Hyperbolic-parabolic
singular perturbations for scalar nonlinearities.  \emph{Appl.\ Anal.}
\textbf{29} (1988), no.~1--2, 19--44.

\bibitem{ghisi1}{\sc M.\ Ghisi}; Global solutions for dissipative
Kirchhoff strings with $m(r)=r^{p}\ (p<1)$.  \emph{J.\ Math.\ Anal.\
Appl.}\ \textbf{250} (2000), no.~1, 86--97.

\bibitem{g:non-lip}{\sc M.\ Ghisi}; Global solutions for dissipative
Kirchhoff strings with non-Lipschitz nonlinear term.  \emph{J.\ Differential
Equations} \textbf{230} (2006), no.~1, 128--139.

\bibitem{gg:k-dissipative}{\sc M.\ Ghisi, M.\ Gobbino}; Global
existence and asymptotic behavior for a mildly degenerate dissipative
hyperbolic equation of Kirchhoff type.  \emph{Asymptot.\  Anal.}\ 
\textbf{40} (2004), no.~1, 25--36.


\bibitem{gg:k-decay}{\sc M.\ Ghisi, M.\ Gobbino}; Hyperbolic-parabolic
singular perturbation for mildly degenerate Kirchhoff equations:
time-decay estimates.  \emph{J.\ Differential Equations} \textbf{245}
(2008), no.~10, 2979--3007.

\bibitem{gg:k-PS}{\sc M.\ Ghisi, M.\ Gobbino}; Hyperbolic-parabolic
singular perturbation for mildly degenerate Kirchhoff equations:
global-in-time error estimates.  To appear on \emph{Commun.\ Pure
Appl.\ Anal.}

\bibitem{gg:cwd}{\sc M.\ Ghisi, M.\ Gobbino}; Hyperbolic-parabolic
singular perturbation for nondegenerate Kirchhoff equations with
critical weak dissipation.  {\tt arXiv:0901.0185  [math.AP]}

\bibitem{k-par}{\sc M.\ Gobbino}; Quasilinear degenerate parabolic
equations of Kirchhoff type.  \emph{Math.\  Methods Appl.\  Sci.}\ 
\textbf{22} (1999), no.~5, 375--388.

\bibitem{k-cattaneo}{\sc M.\ Gobbino}; Singular perturbation
hyperbolic-parabolic for degenerate nonlinear equations of Kirchhoff
type.  \emph{Nonlinear Anal.}\ \textbf{44} (2001), no.~3, Ser.\  A:
Theory Methods, 361--374.

\bibitem{yamazaki} {\sc H.\ Hashimoto, T.\ Yamazaki};
Hyperbolic-parabolic singular perturbation for quasilinear equations
of Kirchhoff type.  \emph{J.\ Differential Equations} \textbf{237}
(2007), no.~2, 491--525.


\bibitem{mizu-ade}{\sc T.\ Mizumachi}; Decay properties of solutions
to degenerate wave equations with dissipative terms.  \emph{Adv.\ 
Differential Equations} \textbf{2} (1997), no.~4, 573--592.

\bibitem{mizu-nc}{\sc T.\ Mizumachi}; Time decay of solutions to
degenerate Kirchhoff type equation.  \emph{Nonlinear Anal.}\ 
\textbf{33} (1998), no.~3, 235--252.

\bibitem{nakao}{\sc M.\ Nakao, J.\ Bae}; Global existence and decay to
the initial-boundary value problem for the Kirchhoff type quasilinear
wave equation with a nonlinear localized dissipation.  \emph{Adv.\
Math.\ Sci.\ Appl.}\ \textbf{13} (2003), no.~1, 165--177.

\bibitem{ny}{\sc K.\ Nishihara, Y.\ Yamada}; On global solutions of
some degenerate quasilinear hyperbolic equations with dissipative
terms.  \emph{Funkcial.\  Ekvac.}\ \textbf{33} (1990), no.~1, 151--159.

\bibitem{ono-kyushu}{\sc K.\ Ono}; Sharp decay estimates of solutions
for mildly degenerate dissipative Kirchhoff equations.  \emph{Kyushu
J.\ Math.}\ \textbf{51} (1997), no.~2, 439--451.

\bibitem{ono-aa}{\sc K.\ Ono}; On global existence and asymptotic
stability of solutions of mildly degenerate dissipative nonlinear wave
equations of Kirchhoff type.  \emph{Asymptot.\  Anal.}\ \textbf{16}
(1998), no.~3-4, 299--314.

\bibitem{ono-wd}{\sc K.\ Ono}; Global solvability for degenerate
Kirchhoff equations with weak dissipation.  \emph{Math.\ Japon.}\
\textbf{50} (1999), no.~3, 409--413.

\bibitem{ono-mm}{\sc K.\ Ono}; Global existence, asymptotic behaviour,
and global non-existence of solutions for damped non-linear wave
equations of Kirchhoff type in the whole space.  \emph{Math.\ Methods
Appl.\ Sci.}\ \textbf{23} (2000), no.~6, 535--560.


\bibitem{yamada}{\sc Y.\ Yamada}; On some quasilinear wave equations with
dissipative terms.  \emph{Nagoya Math.\ J.}\ \textbf{87} (1982),
17--39.


\bibitem{yamazaki-wd}{\sc T.\ Yamazaki}; Hyperbolic-parabolic singular
perturbation for quasilinear equations of Kirchhoff type with weak
dissipation. \emph{Math.\  Methods Appl.\  Sci.}\ In press.

\bibitem{yamazaki-cwd}{\sc T.\ Yamazaki}; Hyperbolic-parabolic
singular perturbation for quasilinear equations of Kirchhoff type with
weak dissipation of critical power.  Preprint. 
\end{thebibliography}
\end{document}